\documentclass[11pt]{article}
\usepackage{latexsym}

\usepackage[tbtags]{amsmath}
\usepackage{epsfig,amstext,amssymb,amsthm,latexsym}
\usepackage{amssymb,color}

\UseRawInputEncoding

\definecolor{c20}{rgb}{0.,0.7,0.}
\definecolor{c30}{rgb}{0.,0.,1.}
\definecolor{c40}{rgb}{1,0.1,0.7}
\definecolor{c50}{rgb}{1,0,0}

\date{}
\makeatletter 
\@addtoreset{equation}{section}
\makeatother 

\begin{document}
\title{On Correlation Coefficients}

\author{ Alexei Stepanov \thanks{\noindent  Education and Research Cluster ``Institute of High Technology",\  Immanuel Kant  Baltic Federal University, A.Nevskogo 14, Kaliningrad, 236041 Russia,  email: alexeistep45@mail.ru}}

\maketitle
\begin{abstract} 
In the present paper, we discuss the Pearson $\rho$,  Spearman $\rho_S$,  Kendall $\tau$ correlation coefficients and their statistical analogues $\rho_n, \rho_{n,S}$ and $\tau_n$. We  propose a new correlation coefficient $r$ and its statistical analogue $r_n$. The coefficient $r$ is based on Kendal's and Spearman's correlation coefficients. A new extension of the Pearson correlation coefficient  is also discussed. We conduct  simulation experiments and study the behavior of the above correlation coefficients. We  observe  that the behavior of $\rho_n$ can be very different from the behavior of the rank correlation coefficients $\rho_{n,S}, \tau_n$ and $r_n$, which, in  turn, behave in a similar way. The question arises:   which  correlation coefficient  better measures the dependence rate? We try to answer this question in the final conclusion.
\end{abstract}



\section{Introduction} Let  $(X,Y), (X_1,Y_1),\ldots,(X_n,Y_n)$  be independent and identically distributed  random vectors with absolutely continuous bivariate distribution function $F(x,y)$,   density function  $f(x,y)$ and marginal distribution functions    $H(x)$ and $G(y)$. A rate of dependence between  $X$ and $Y$ indicates how closely these variables  are related. This ranges from complete dependence to independence. There are three basic  measures of dependence: the Pearson $\rho$, the Spearman $\rho_S$ and the Kendall $\tau$ correlation coefficients. 

The rate of dependence between the variables $X$ and $Y$ is most often measured by the Pearson correlation coefficient
$$
\rho=\frac{E(X-EX)(Y-EY)}{\sigma_X\sigma_Y}.
$$
The random variable
$$
\rho_n =\frac{\sum\limits_{i=1}^n (X_i-\bar{X})(Y_i-\bar{Y})} {\sqrt{\sum\limits_{i=1}^n (X_i-\bar{X})^2 \sum\limits_{i=1}^n (Y_i-\bar{Y})^2}}
$$
is known as the sample Pearson correlation coefficient. It is a good approximation for  $\rho $ for large  $n$ since $\rho_n\stackrel{p}{\rightarrow}\rho$; see, for example,  Fisher (1921). 

Let  $X_{(1)}\leq\ldots\leq X_{(n)}$ be the order statistics obtained from the sample   $X_1,\ldots,X_n$. For these order statistics, let us define their concomitants   $Y_{[1]},\ldots,Y_{[n]}$. Let $X_i=X_{(j)}$, then  $Y_{[j]}=Y_i$ be the concomitant of the order statistic  $X_{(j)}$.  The concept of concomitants was  proposed and discussed in the papers of   David and Galambos (1974), Bhattacharya (1974, 1984), Egorov and Nevzorov (1984),  Goel and Hall (1994), David (1994), Chu et al. (1999), David and Nagaraja (2003), Balakrishnan and Lai (2009), Bairamov and Stepanov (2010),   Balakrishnan and Stepanov (2015). See also the references therein. Let  
$$
R_j=\sum_{i=1}^{j-1}I(Y_{[i]}\leq Y_{[j]})\quad (1\leq i<j\leq n),
$$
where $I$ is the indicator function, be  the rank of $Y_{[j]}$  amongst the concomitants  $Y_{[1]},\ldots, Y_{[j]}$. The random variable  
$$
\tau _n=\frac{4\sum_{j=2}^nR_{j}}{n(n-1)}-1
$$
is known as the rank Kendall correlation coefficient. The theoretical analogue of $\tau_n$ is presented by
$$
\tau=4E[F(X,Y)]-1=4\int_{\mathbb{R}^2}F(x,y)f(x,y)dxdy-1.
$$
It is known, see, for example, Balakrishnan and Lai (2009) or Stepanov (2024), that $\tau=E\tau_n$ and $\tau_n\stackrel{p}{\rightarrow}\tau$. Let 
$$
R_{j,n}=\sum_{i=1}^{n}I(Y_{[i]}\leq Y_{[j]})\quad (1\leq j\leq n)
$$
be  the rank of $Y_{[j]}$  amongst the concomitants  $Y_{[1]},\ldots, Y_{[n]}$. The random variable  
$$
\rho_{n,S} =1-\frac{6\sum_{j=1}^n(R_{j,n}-j)^2}{n^3-n}
$$ 
is known as the rank Spearman correlation coefficient. The theoretical analogue of $\rho_{n,S}$ is presented by
$$
\rho_S =12E[H(X)G(Y)]-3=12\int_{\mathbb{R}^2}H(x)G(y)f(x,y)dxdy-3.
$$ 
The random variable $F(X,Y)-H(X)G(Y)$ is called the measure of quadrant dependence; see, for example, p. 158 of Balakrishnan and Lai (2009).  

The introduced above correlation coefficients $\rho,\ \rho_S$ and $\tau$ as well as their statistical analogues $\rho_n,\ \rho_{n,S}$ and $\tau_n$ have some advantages and disadvantages. Which correlation coefficient better measures the dependence rate? For example, it is known that in many cases $\mid \tau\mid\leq \mid \rho_S \mid \leq \mid \rho\mid$ and $\mid \tau_n\mid\leq \mid \rho_{n,S} \mid \leq \mid \rho_n\mid$. Does it mean that, say,  $\rho$ is better  than $\tau$? It is also known, for example, that $\rho_n$ is very sensitive to contamination of outliers, while $\rho_{n,S}$ and $\tau_n$ are rather robust. This property  can be viewed as an obvious disadvantage of $\rho_n$ and, correspondingly, of $\rho$. Another disadvantage of $\rho$ and, correspondingly of $\rho_n$, is that $\rho$ cannot be defined  if  the  second moments do not exist. Observe that $\rho_S$ and $\tau$ exist for any absolutely continuous non-degenerate distribution function. Another possible disadvantage of $\rho$  is that it is defined only in terms of the first and second moments and, accordingly, loses much information about the rate of dependence  between  $X$ and $Y$. On the other hand, this can be viewed as an advantage since one can find $\rho$ knowing only the first and second moments and does not need to know $F$. Since $\tau_n$ and $\rho_{n,S}$ measure the dependence rate between ranks, rather than between actual values of $X$ and $Y$, the  correlation coefficients $\tau$ and $\rho_{S}$ are unaffected by any increasing transformation of $X$ and $Y$, for example, by the monotone transformation based on  the marginals of $F$. Observe that $\rho$ is not invariant under this transformation. Then, if a distribution function $F$ depends on the parameter $t$, and if this parameter can be suppressed by one of the marginal transformations $u=H(x)$ or $v=G(y)$, then  $\rho_{S}$ and $\tau$ are free of $t$ for this distribution function $F$, wherein $\rho$ can depend on $t$. In that case, the coefficients $\rho_S$ and $\tau$ are not elastic and, it seems that  $\rho(t)$, which is elastic here,  measures the dependence rate better; see Example~3.3 below. 

Although it is customary in bivariate data analysis to compute the dependence rate by using only one correlation coefficient, in this paper we propose to use a new correlation coefficient, which is a linear combination of $\rho_S$ and $\tau$. In general, before applying any correlation coefficient, a researcher should study the scatter plot of data and decide which  correlation coefficient can be most helpful in that situation. The coefficient $\rho_n$ can be better used when the scatter plot indicates a linear relationship between $X$ and $Y$, and $\rho_{n,S}$ and $\tau_n$ can be  exploited in  other cases.  

The rest of our paper is organized as follows. In Section~2, we propose a new correlation coefficient, which is based on Kendal's and Spearman's correlation coefficients. There,  we also give a new extension of the Pearson correlation coefficient. Section~3 contains examples, which illustrate the work of all four correlation coefficients that we discuss in our paper. In the end of this paper we  present a conclusion.

\section{New Correlation Coefficient. Extension of  Pearson Correlation Coefficient.}
As was already pointed out, the following question is interesting: which  correlation coefficient measures the dependence rate better?  Before answering it, we analyze the known coefficients $\rho,\ \rho_S$ and $\tau$ and introduce a new one $r$.

Let us consider the dependence function $F(x,y)-H(x)G(y)$, which gives us the dependence value in every point $(x,y)$ of the distribution support. We adjust the definition of the dependence function  and redefine it as $D(x,y)=6(F(x,y)-H(x)G(y))$. We then can view any correlation coefficient     as  some average values of  the dependence rate, or of the dependence function $D$. Observe that the dependence function does not take into account the likelihood of obtaining points $(x,y)$. Motivated by the above argument, we assume that $r=E(D(X,Y))$ represents  the ``ideal" average of the dependence function. 

\noindent {\bf New Correlation Coefficient}\ In this paper, we  propose to consider 
\begin{eqnarray*}
r&=&E(D(X,Y))\\
&=& 6E[F(X,Y)-H(X)G(Y)]=6\int_{\mathbb{R}^2}(F(x,y)-H(x)G(y))f(x,y)dxdy
\end{eqnarray*}
as a new theoretical correlation coefficient. Observe that $r=\frac{3\tau-\rho_S}{2}$. The following two inequalities 
$$
\frac{3\tau-1}{2}\leq \rho_S\ (\mbox{if}\ \tau\geq0),\quad   \rho_S\leq \frac{1+3\tau}{2}\ (\mbox{if}\ \tau\leq 0)
$$
are given in (4.18) of Balakrishnan and Lai (2009). The first inequality implies that $r\leq 1$, when the second one assumes that $r\geq -1$. Since in a statistical experiment we usually do not know the underlying $F$, we need to propose a statistical analogue of the  correlation coefficient $r$.  Taking into account the definition of $r$, we propose to view
$$
r_n=1.5 \tau_n-0.5 \rho_{n,S}
$$
as the corresponding rank correlation coefficient. It is possible that $r_n$ can be presented in a more simple form. It was mentioned above  that  $E\tau_n=\tau$. On page 157 of Balakrishnan and Lai (2009) one can find that $E\rho_{n,S}=\frac{(n-2)\rho_S+3\tau}{n+1}$.  It follows that 
$$
Er_n=\frac{3\tau}{2}-\frac{(n-2)\rho_S+3\tau}{2(n+1)}=\frac{3n\tau-(n-2)\rho_S}{2(n+1)}\not=1.5\tau-0.5\rho_S=r.
$$
That way, although $r_n$ is a biased estimator of $r$, it is asymptotically unbiased. The central limit theorem for $\rho_{n,S}$ was proved in Cifarelli et al. (1996), from which we gain that $\rho_{n,S}\stackrel{p}{\rightarrow}\rho_S$. It is also known, see, for example, Stepanov (2024), that $\tau_{n}\stackrel{p}{\rightarrow}\tau\ (n\rightarrow \infty)$. By Slutsky's theorem, we obtain that $r_n\stackrel{p}{\rightarrow}r$. We see that $r_n$ can be considered as a statistical analogue of $r$. 

Since $r$ is a mixture of $\rho_S$ and $\tau$, it inherits all advantages and disadvantages of them. It is rather robust to the influence of outliers, but it is constant  when $F$ depends on a  parameter $t$ and when $t$ can be suppressed by the marginal transformations. Another situation when $r$ presumably works unsatisfactory (in comparison with $\rho$) is when the dependence rate is close to $\pm 1$.  By the definition of $r$, it is fair to compare $r_n$ only with   $\tau_n$ and $\rho_{n,S}$, because these three coefficients   measure  dependence  between the ranks. As was pointed out, $r$ represents some kind of the ``ideal" average of the dependence function $D$. By this argument, $r$ is better than $\rho_S$ and $\tau$. The  correlation coefficient $\rho$ represents the average of the random variable $\frac{XY-E(XY)}{\sigma_X\sigma_Y}$.  Accordingly, the coefficient $\rho$ is simpler, but some information about dependence can be lost, and $\rho$ cannot be  so  accurate when the situation is ``far" from the linear relationship between $X$ and $Y$. The coefficient $\rho_n$ measures  dependence between sample members and is somewhat different from the rank correlation coefficients $\rho_{n,S},\ \tau_n$ and $r_n$, which measure association between the ranks.

An interesting fact should be highlighted. By analyzing the dependence function $D$, we conclude that the variables $X$ and $Y$ are independent iff $Var(D(X,Y))=0$. We do not know how this discovery can be further used.

It should be noted that correlation coefficients can be compared by the sample variances obtained by generation of these correlation coefficients. If a sample variance is low, it indirectly  indicates  that the corresponding correlation coefficient  measures the dependence rate accurately, i.e., with  small   squared error.

Let $\bar{\rho}_n$ and $S^2_{\rho_n}$ denote the  mean and  variance of a sample consisting of generated Pearson correlation coefficients. Let the designations for the other correlation coefficients  are similar. Looking at the results of simulation in Section~3, we find  that generally 
$$
|\bar{r}_n|\leq |\bar{\tau}_n | \leq | \bar{\rho}_{n,S} |\leq |\bar{\rho}_n |,
$$ 
and the inequalities 
$$
S^2_{r_{n}} \leq S^2_{\tau_{n}}\leq S^2_{\rho_{n,S}}\leq S^2_{\rho_{n}}
$$
hold true with the exception for the bivariate normal distribution, when $t$ is close to $\pm 1$.  As we can see  in Section~3 below, the rank correlation coefficient $r_n$ inherits all the properties of $\rho_S$ and $\tau$, but has one advantage over them. Most often  it has the smallest sample variance. Comparing the results of simulation in Section~3, we find that $\rho_{n,S},\ \tau_n$, and now $r_n$,  ``underestimate"\ the dependence rate in comparison with $\rho_n$. 

In the end of this section, we would like to propose a new extension of the Pearson correlation coefficient $\rho$. In the case when $Var (X)$ (or $Var (Y)$) does not exist, the coefficient $\rho$  is formally not defined. It was known that  it can be generalized and redefined, for example, for the bivariate stable distribution functions; see p. 146 of Balakrishnan and Lai (2007). Thus, for the symmetric bivariate  Cauchy distribution function, we  can apply the Cauchy principal value. Then, for this distribution function,  $cov(X,Y)=0$ and we can put $\rho=0$.

In this work, we propose a new extension of $\rho$. In the case, when the second moments do not exist, but simulation experiments show that $\bar{\rho}_n\approx \tilde{\rho}$ for large $n$, we propose to consider $\tilde{\rho}$ as a generalization of $\rho$ provided that the values of $\bar{\rho}_n$ are ``close" to the values of  $\bar{\rho}_{n,S},\ \bar{\tau}_n$ and $\bar{r}_n$; see Examples 3.2 below.

\section{Examples}

\noindent {\bf Example~3.1}\ Let
$$
F_t(x,y)=\frac{1}{2\pi\sqrt{1-t^2}}\int_{-\infty}^x\int_{-\infty}^y e^{-\frac{u^2-2t uv+v^2}{2(1-t^2)}}dudv\quad (x,y\in \mathbb{R},\ -1<t<1)
$$
be a bivariate normal distribution. We have 
$$
\rho(t)=t,\quad \rho_S(t)=\frac{6}{\pi}asin(t/2),\quad \tau(t)=\frac{2}{\pi}asin(t)\quad r(t)=\frac{3}{\pi}\left(asin(t)-asin(t/2)\right).
$$ 
For visual comparison, we present the above correlation coefficients for different values of  $t$  in Table~3.1.  

\begin{table}[ht]
\begin{center}
{\em Table~3.1 Values of correlation coefficients for different $t$.}

\vspace*{1ex}
\begin{tabular}{|c|c|c|c|c|c|c|c|c|}
\hline
&  $t=-0.99$ & $t=-0.70$ & $t=-0.30$ & $t=-0.10$ & $t=0.10$ & $t=0.30$ & $t=0.70$ & $t=0.99$\\
\hline
$\rho$  & -0.9900 & -0.7000 & -0.3000 & -0.1000 & 0.1000 & 0.3000 & 0.7000 & 0.9900\\
\hline
$\rho_s$& -0.9890 & -0.6829 & -0.2876 & -0.0955 & 0.0955 & 0.2876 & 0.6829 & 0.9890\\
\hline
$\tau$  & -0.9098 & -0.4936 & -0.1940 &  -0.0638 & 0.0638 & 0.1940 & 0.4936 & 0.9098 \\
\hline
$r$     & -0.8703 & -0.3990 & -0.1472 & -0.0479 & 0.0479 & 0.1472 & 0.3990 & 0.8703\\
\hline
\end{tabular}
\end{center}
\end{table}

It is clear from Table~3.1 that $\mid r |\leq | \tau | \leq  | \rho_S | \leq | \rho\mid$. We  conducted several simulation experiments.  We   generated one thousand times these correlation coefficients for samples $X_1,\ldots, X_{1000}$ and $Y_1,\ldots, Y_{1000}$, computed their means  and variances and presented them in  Tables~3.2--3.3, respectively.

\begin{table}[ht]
\begin{center}
{\em Table~3.2 Sample means of correlation coefficients for different $t$.}

\vspace*{1ex}
\begin{tabular}{|c|c|c|c|c|c|c|c|c|}
\hline
&  $t=-0.99$ & $t=-0.70$ & $t=-0.30$ & $t=-0.10$ & $t=0.10$ & $t=0.30$ & $t=0.70$ & $t=0.99$\\
\hline
$\bar{\rho}_{n}$   & -0.9900 & -0.7000 & -0.3002 & -0.0999 & 0.0996 & 0.3003 & 0.6985 & 0.9900\\
\hline
$\bar{\rho}_{n,S}$ & -0.9887 & -0.6825 & -0.2874 & -0.0955 & 0.0949 & 0.2878 & 0.6812 & 0.9888\\
\hline
$\bar{\tau}_{n}$   & -0.9100 & -0.4940 & -0.1940 & -0.0638 & 0.0634 & 0.1944 & 0.4926 & 0.9100\\
\hline
$\bar{r}_{n}$      & -0.8706 & -0.3994 & -0.1473 & -0.0480 & 0.0477 & 0.1476 & 0.3982 & 0.8706\\
\hline
\end{tabular}
\end{center}
\end{table}

\begin{table}[ht]
\begin{center}
{\em Table~3.3 Sample variances of correlation coefficients for different $t$.}

\vspace*{1ex}
\begin{tabular}{|c|c|c|c|c|c|c|c|c|}
\hline
                   &  $t=-0.99$         & $t=-0.70$          & $t=-0.30$          & $t=-0.10$          & $t=0.10$           & $t=0.30$           & $t=0.70$           & $t=0.99$\\
\hline
$S^2_{\rho_{n}}$   & $4.1\cdot 10^{-7}$ & $2.9\cdot 10^{-4}$ & $8.1\cdot 10^{-4}$ & $1.0\cdot 10^{-3}$ & $9.0\cdot 10^{-4}$ & $8.4\cdot 10^{-4}$ & $2.6\cdot 10^{-4}$ & $4.0\cdot 10^{-7}$\\
\hline
$S^2_{\rho_{n,S}}$ & $7.0\cdot 10^{-7}$ & $3.6\cdot 10^{-4}$ & $8.5\cdot 10^{-4}$ & $1.0\cdot 10^{-3}$ & $9.0\cdot 10^{-4}$ & $8.8\cdot 10^{-4}$ & $3.2\cdot 10^{-4}$ & $7.2\cdot 10^{-7}$ \\
\hline
$S^2_{\tau_n}$     & $9.9\cdot 10^{-6}$ & $2.6\cdot 10^{-4}$ & $4.0\cdot 10^{-4}$ & $4.5\cdot 10^{-4}$ & $4.1\cdot 10^{-4}$ & $4.2\cdot 10^{-4}$ & $2.3\cdot 10^{-4}$ & $9.8\cdot 10^{-6}$ \\
\hline
$S^2_{r_n}$        & $1.8\cdot 10^{-5}$ & $2.2\cdot 10^{-4}$ & $2.4\cdot 10^{-4}$ & $2.5\cdot 10^{-4}$ & $2.3\cdot 10^{-4}$ & $2.6\cdot 10^{-4}$ & $1.9\cdot 10^{-4}$ & $1.8\cdot 10^{-5}$ \\
\hline
\end{tabular}
\end{center}
\end{table}
It seems that $r_n$ is a good rank correlation coefficient, since it has the uniformly small  sample variance.   

We generated  once more one thousand times correlation coefficients $\rho_n,\ \rho_{n,S},\ \tau_n$ and $r_n$.  One thousand observations were taken from  $F$ and five observations were outliers. We obtained the following results for $t=-0.99$
$$
\bar{\rho}_n=-0.3305,\ \bar{\rho}_{n,S}=-0.9592,\ \bar{\tau}_n=-0.8910,\ \bar{r}_n=-0.8569.
$$
As we can see, the effect of outliers on $\rho_n$  is great. 

We also simulated the correlation coefficients that measure correlation between the correlation coefficients themselves. Let, for example, $\rho_{100}^{(i)}, \rho_{100,S}^{(j)}\ (1\leq i,j\leq 100)$ be samples of size 100 of Pearson's and Spearman's correlation coefficients and $\bar{\rho}_{100}(\rho_{100}^{(i)}, \rho_{100,S}^{(j)})$ be the mean of Pearson's correlation coefficients, which measure the dependence rate between the samples consisting of Pearson's and Spearman's correlation coefficients. For $t=-0.99$, we present  all  results for all types of correlation coefficients:
$$
\bar{\rho}_{100}(\rho_{100}^{(i)}, \rho_{100,S}^{(j)})=0.7340,\quad \bar{\rho}_{100}(\rho_{100}^{(i)}, \tau_{100}^{(j)})=0.7987,\quad \bar{\rho}_{100}(\rho_{100}^{(i)}, r_{100}^{(j)})=0.8017,
$$
$$
\bar{\rho}_{100}(\rho_{100,S}^{(i)}, \tau_{100}^{(j)})=0.9578,\quad \bar{\rho}_{100}(\rho_{100,S}^{(i)}, r_{100}^{(j)})=0.9485,\quad \bar{\rho}_{100}(\tau_{100}^{(i)}, r_{100}^{(j)})=0.9985,
$$
$$
\bar{\rho}_{100,S}(\rho_{100}^{(i)}, \rho_{100,S}^{(j)})=0.7242,\quad \bar{\rho}_{100,S}(\rho_{100}^{(i)}, \tau_{100}^{(j)})=0.7826,\quad \bar{\rho}_{100,S}(\rho_{100}^{(i)}, r_{100}^{(j)})=0.7856,
$$
$$
\bar{\rho}_{100,S}(\rho_{100,S}^{(i)}, \tau_{100}^{(j)})=0.9568,\quad \bar{\rho}_{100,S}(\rho_{100,S}^{(i)}, r_{100}^{(j)})=0.9479,\quad \bar{\rho}_{100,S}(\tau_{100}^{(i)}, r_{100}^{(j)})=0.9992,
$$
$$
\bar{\tau}_{100}(\rho_{100}^{(i)}, \rho_{100,S}^{(j)})=0.5385,\quad \bar{\tau}_{100}(\rho_{100}^{(i)}, \tau_{100}^{(j)})=0.5972,\quad \bar{\tau}_{100}(\rho_{100}^{(i)}, r_{100}^{(j)})=0.5973,
$$
$$
\bar{\tau}_{100}(\rho_{100,S}^{(i)}, \tau_{100}^{(j)})=0.9578,\quad \bar{\tau}_{100}(\rho_{100,S}^{(i)}, r_{100}^{(j)})=0.8121,\quad \bar{\tau}_{100}(\tau_{100}^{(i)}, r_{100}^{(j)})=0.9861.
$$
As we can see the correlations between the rank correlation coefficients are rather close to one. The Pearson sample correlation coefficient is not so much correlated with the rank correlation coefficients. The smallest and largest  correlation coefficients  for the normal case for $t=-0.99$ are $\bar{\tau}_{100}(\rho_{100}^{(i)}, \rho_{100,S}^{(j)})=0.5385$ and $\bar{\rho}_{100,S}(\tau_{100}^{(i)}, r_{100}^{(j)})=0.9992$, respectively. The smallest and largest  correlation coefficients   for $t=-0.7$ are, correspondingly,  $\bar{\tau}_{100}(\rho_{100}^{(i)}, \rho_{100,S}^{(j)})=0.7126$ and $\bar{\rho}_{100,S}(\tau_{100}^{(i)}, r_{100}^{(j)})=0.9956$. As $t$ advances to zero the difference between the largest and smallest   correlation coefficients  decreases. Thus, for $t=-0.1$, we have  $\bar{\tau}_{100}(\rho_{100}^{(i)}, r_{100}^{(j)})=0.7786$ and $\bar{\rho}_{100,S}(\rho_{100,S}^{(i)}, \tau_{100}^{(j)})=0.9956$. For positive $t$ the situation is symmetric. We  conclude that there is  some similarity in behaviors of $\rho_n$ and the rank correlation coefficients when $t$ is not close to $\pm 1$.

\noindent {\bf Example~3.2}\ Let us consider the bivariate Pareto distribution function
$$
F_t(x,y)=1-\frac{1}{(y+1)^t}-\frac{1}{(x+1)^t}+\frac{1}{(x+y+1)^t}\quad (x>0,\ y>0,\ t>0).
$$
We have $\rho(t) =\frac{1}{t}\ (t>2)$ and   $\tau(t)=\frac{1}{2t+1}\ (t>0)$. Since  correlation coefficients $\rho_s(t)$ and $r(t)$ can be found  only numerically, we present   in Table~3.4   their values for different  $t$ along with the values of $\rho(t)$ and   $\tau(t)$. The means  and variances of these correlation coefficients are  presented in  Tables~3.5--3.6, respectively.
\begin{table}[ht]
\begin{center}
{\em Table~3.4 Values of correlation coefficients for different $t$.}

\vspace*{1ex}
\begin{tabular}{|c|c|c|c|c|c|c|c|c|c}
\hline
  &  $t=0.1$ & $t=0.5$ & $t=1$ & $t=2.1$ & $t=5$ & $t=10$ & $t=50$ & $t=100$\\
 \hline
 $\rho$  & -- & -- & -- & 0.4761 & 0.2000& 0.1000 & 0.0200 & 0.0100\\
\hline
$\rho_s$ & 0.92368 & 0.6822 & 0.4784 & 0.2839 & 0.1358 & 0.0714 & 0.0149 & 0.0075\\
\hline
$\tau$  & 0.8333 & 0.5000 & 0.3333 &  0.1923 & 0.0909 & 0.0476 & 0.0099 & 0.0050\\
\hline
$r$    & 0.6903 & 0.4089 & 0.2608 & 0.1465 & 0.0684 & 0.0357 & 0.0074 & 0.0037\\
\hline
\end{tabular}
\end{center}
\end{table}

\begin{table}[ht]
\begin{center}
{\em Table~3.5 Sample means  of correlation coefficients for different $t$.}

\vspace*{1ex}
\begin{tabular}{|c|c|c|c|c|c|c|c|c|c}
\hline
 & $t=0.1$ & $t=0.5$ & $t=1$ & $t=2.1$ & $t=5$ & $t=10$ & $t=50$ & $t=100$\\
 \hline
 $\bar{\rho}_{n}$  & 0.9545 & 0.7972 & 0.6502 & 0.4401 & 0.1996 & 0.1003 & 0.0198 & 0.0099\\
\hline
$\bar{\rho}_{n,S}$  & 0.9580 & 0.6818 & 0.4784 & 0.2831 & 0.1364 & 0.0709 & 0.0143 & 0.0073\\
\hline
$\bar{\tau}_{n}$   & 0.8334 & 0.5000 & 0.3337 & 0.1919 & 0.0914 & 0.0474 & 0.0096 & 0.0049\\
\hline
 $\bar{r}_{n}$    & 0.7712 & 0.4091 & 0.2613 & 0.1463 & 0.0689 & 0.0356 & 0.0072 & 0.0037\\
\hline
\end{tabular}
\end{center}
\end{table}
\begin{table}[ht]
\begin{center}
{\em Table~3.6 Sample variances  of correlation coefficients for different $t$.}

\vspace*{1ex}
\begin{tabular}{|c|c|c|c|c|c|c|c|c|c}
\hline
  &  $t=0.1$ & $t=0.5$ & $t=1$ & $t=2.1$ & $t=5$ & $t=10$ & $t=50$ & $t=100$\\
 \hline
 $S^2_{\rho_{n}}$  &$1.9\cdot 10^{-2}$ & $5.4\cdot 10^{-2}$ & $5.2\cdot 10^{-2}$ & $2.2\cdot 10^{-2}$ & $3.5\cdot 10^{-3}$ & $1.6\cdot 10^{-3}$ & $9.5\cdot 10^{-4}$ & $1.0\cdot 10^{-3}$ \\
\hline
$S^2_{\rho_{n,S}}$  & $1.4\cdot 10^{-5}$ & $4.0\cdot 10^{-4}$ & $7.1\cdot 10^{-4}$   & $9.6\cdot 10^{-4}$ & $11.0\cdot 10^{-4}$ & $9.7\cdot 10^{-4}$ & $9.9\cdot 10^{-4}$ & $1.0\cdot 10^{-3}$ \\
\hline
$S^2_{\tau_n}$    & $4.9\cdot 10^{-5}$ & $2.9\cdot 10^{-4}$ & $3.9\cdot 10^{-4}$ & $4.7\cdot 10^{-4}$ & $5.0\cdot 10^{-4}$ & $4.3\cdot 10^{-4}$ & $4.4\cdot 10^{-4}$ & $4.5\cdot 10^{-4}$ \\
\hline
 $S^2_{r_n}$  & $7.6\cdot 10^{-5}$ & $2.5\cdot 10^{-4}$ & $2.7\cdot 10^{-4}$ & $2.8\cdot 10^{-4}$ & $2.9\cdot 10^{-4}$ & $2.5\cdot 10^{-4}$ & $2.5\cdot 10^{-4}$ & $2.5\cdot 10^{-4}$ \\
\hline
\end{tabular}
\end{center}
\end{table}
It seems that again  $r$  has the smallest sample variance and can be viewed as a solid  correlation coefficient.  

It should be noted  that $\rho_n$ measures the dependence rate when $t\in(0,2]$ and when its theoretical analogue $\rho$ does not exist. Moreover, it follows the trend given by the rank correlation coefficients, i.e., $\bar{\rho}_n(t)$ tends to one as $t$ tends to zero. Taking into account the last paragraph of introduction, we define the correlation coefficient $\tilde{\rho}$, which is an extension of  $\rho$. Let $\tilde{\rho}$ be such that  $\tilde{\rho}(t)=\rho(t)$ if $t>2$ and  $\tilde{\rho}(t)= \bar{\rho}_n$, where $\bar{\rho}_n(t)$ is found  by simulation for every $t\in(0,2]$.

We also simulated the correlation coefficients which measure correlation between the correlation coefficients themselves. For $t=0.1$ we present below all  results for all  types of correlation coefficients:
$$
\bar{\rho}_{100}(\rho_{100}^{(i)}, \rho_{100,S}^{(j)})=-2.8\cdot 10^{-3},\quad \bar{\rho}_{100}(\rho_{100}^{(i)}, \tau_{100}^{(j)})=-3.9\cdot 10^{-4},\quad \bar{\rho}_{100}(\rho_{100}^{(i)}, r_{100}^{(j)})=1.5\cdot 10^{-4},
$$
$$
\bar{\rho}_{100}(\rho_{100,S}^{(i)}, \tau_{100}^{(j)})=0.9655,\quad \bar{\rho}_{100}(\rho_{100,S}^{(i)}, r_{100}^{(j)})=0.9489,\quad \bar{\rho}_{100}(\tau_{100}^{(i)}, r_{100}^{(j)})=0.9983,
$$
$$
\bar{\rho}_{100,S}(\rho_{100}^{(i)}, \rho_{100,S}^{(j)})=-4.2\cdot 10^{-3},\quad \bar{\rho}_{100,S}(\rho_{100}^{(i)}, \tau_{100}^{(j)})=-2.4\cdot 10^{-3},\quad \bar{\rho}_{100,S}(\rho_{100}^{(i)}, r_{100}^{(j)})=-2.8\cdot 10^{-3},
$$
$$
\bar{\rho}_{100,S}(\rho_{100,S}^{(i)}, \tau_{100}^{(j)})=0.9646,\quad \bar{\rho}_{100,S}(\rho_{100,S}^{(i)}, r_{100}^{(j)})=0.9491,\quad \bar{\rho}_{100,S}(\tau_{100}^{(i)}, r_{100}^{(j)})=0.9976,
$$
$$
\bar{\tau}_{100}(\rho_{100}^{(i)}, \rho_{100,S}^{(j)})=-2.7\cdot 10^{-3},\quad \bar{\tau}_{100}(\rho_{100}^{(i)}, \tau_{100}^{(j)})=-1.4\cdot 10^{-3},\quad \bar{\tau}_{100}(\rho_{100}^{(i)}, r_{100}^{(j)})=-1.6\cdot 10^{-3},
$$
$$
\bar{\tau}_{100}(\rho_{100,S}^{(i)}, \tau_{100}^{(j)})=0.9655,\quad \bar{\tau}_{100}(\rho_{100,S}^{(i)}, r_{100}^{(j)})=0.8132,\quad \bar{\tau}_{100}(\tau_{100}^{(i)}, r_{100}^{(j)})=0.9693.
$$
As we can see again, the correlations between the rank correlation coefficients are rather close to one. On the contrary,  $\rho_n$  is almost uncorrelated with the rank correlation coefficients, though the value of $\bar{\rho}_{n}$ is rather close to the values of $\bar{\rho}_{n,S},\ \bar{\tau}_{n}$ and $\bar{r}_{n}$. When $t$ increases the situation becomes more even. For $t=1, 10, 100, 1000$  we present  only the smallest and  largest  correlation coefficients:
$$
t=1,\quad \bar{\tau}_{100}(\rho_{100}^{(i)}, \rho_{100,S}^{(j)})=0.0358,\quad \bar{\tau}_{100}(\tau_{100}^{(i)}, r_{100}^{(j)})=0.9961,
$$
$$
t=10,\quad \bar{\tau}_{100}(\rho_{100}^{(i)}, r_{100,S}^{(j)})=0.4558,\quad \bar{\rho}_{100}(\tau_{100}^{(i)}, r_{100}^{(j)})=0.9977,
$$
$$
t=100,\quad \bar{\tau}_{100}(\rho_{100}^{(i)}, r_{100,S}^{(j)})=0.5507,\quad \bar{\rho}_{100}(\tau_{100}^{(i)}, r_{100}^{(j)})=0.9977.
$$
$$
t=1000,\quad \bar{\tau}_{100}(\rho_{100}^{(i)}, r_{100,S}^{(j)})=0.5607,\quad \bar{\rho}_{100}(\tau_{100}^{(i)}, r_{100}^{(j)})=0.9979.
$$
As $t$  increases further, the difference between the largest and  smallest  correlation coefficients does not change much.

\noindent {\bf Example~3.3}\ Let
$$
F_t(x,y)=1-e^{-x}-\frac{1}{(1+y)^t}+\frac{e^{-x(1+y)^t}}{(1+y)^t}\quad (x,y>0).
$$
We have   $\rho =-\frac{\sqrt{t(t-2)}}{2t-1}\ (t>2)$. For  $t>0$ the rank correlation coefficients are free of $t$ and take the following values
$$
\rho_S=12\log 2-9=-0.6822,\quad  \tau=-0.5,\quad  r= 15/4-6\log 2=-0.4089.
$$
One can see that here only $\rho$ depends on $t$ and presumably  better reflects the dependence rate when $t>2$.  The values of $\rho$ along with the means of all correlation coefficients are presented  in Table~3.7  for $t>0$. The corresponding sample variances are presented in Table~3.8.

\begin{table}[ht]
\begin{center}
{\em Table~3.7 Values of $\rho$ and sample means of  correlation coefficients.}

\vspace*{1ex}
\begin{tabular}{|c|c|c|c|c|c|c|c|c|}
\hline
 &  $t=0.1$ & $t=0.5$ & $t=1$ & $t=2.1$ & $t=5$ & $t=10$ & $t=50$ & $t=100$\\
 \hline
 $\rho$            & --      & --       &  --     & -0.1432 & -0.4303 & -0.4708 & -0.4948 & -0.4975\\
 \hline
 $\bar{\rho}_{n}$  & -0.0324 & -0.0525  & -0.1192 & -0.2998 & -0.4364 & -0.4726 & -0.4961 & -0.4989\\
\hline
$\bar{\rho}_{n,S}$ & -0.6820 & -0.6816 & -0.6814 & -0.6819 & -0.6815 & -0.6811 & -0.6814 & -0.6820\\
\hline
$\bar{\tau}_{n}$   & -0.5002 & -0.4999 & -0.4998 & -5.002 & -0.4999 & -0.4995 & -0.4998 & -0.5002\\
\hline
 $\bar{r}_{n}$     & -0.4093 & -0.4091 & -0.4090 & -0.4093 & -0.4091& -0.4087 & -0.4089 & -0.4093\\
\hline
\end{tabular}
\end{center}
\end{table}

\begin{table}[ht]
\begin{center}
{\em Table~3.8 Sample variances  of correlation coefficients for different $t$.}

\vspace*{1ex}
\begin{tabular}{|c|c|c|c|c|c|c|c|c|}
\hline
 &  $t=0.1$ & $t=0.5$ & $t=1$ & $t=2.1$ & $t=5$ & $t=10$ & $t=50$ & $t=100$\\
 \hline
 $S^2_{\rho_{n}}$  & $2.5 \cdot 10^{-5}$ & $3.7\cdot 10^{-4}$ & $2.7\cdot 10^{-3}$ & $4.2\cdot 10^{-3}$ & $7.8\cdot 10^{-4}$ & $3.7\cdot 10^{-4}$ & $2.5\cdot 10^{-4}$ & $2.5\cdot 10^{-4}$\\
\hline
$S^2_{\rho_{n,S}}$ & $3.7\cdot 10^{-4}$ & $3.7\cdot 10^{-4}$ & $3.8\cdot 10^{-4}$ & $3.8\cdot 10^{-4}$ & $3.8\cdot 10^{-4}$ & $3.8\cdot 10^{-4}$ & $3.7\cdot 10^{-4}$ & $3.8\cdot 10^{-4}$\\
\hline
$S^2_{\tau_n}$     & $2.7\cdot 10^{-4}$ & $2.7\cdot 10^{-4}$ & $2.7\cdot 10^{-4}$ & $2.7\cdot 10^{-4}$ & $2.7\cdot 10^{-4}$ & $2.7\cdot 10^{-4}$ & $2.6\cdot 10^{-4}$ & $2.7\cdot 10^{-4}$\\
\hline
  $S^2_{r_n}$      & $2.3\cdot 10^{-4}$ & $2.3\cdot 10^{-4}$ & $2.3\cdot 10^{-4}$ & $2.3\cdot 10^{-4}$ & $2.3\cdot 10^{-4}$ & $2.3\cdot 10^{-4}$ & $2.2\cdot 10^{-4}$ & $2.3\cdot 10^{-4}$\\
\hline
\end{tabular}
\end{center}
\end{table} 
The coefficient $\rho_n$ again measures the dependence rate for $t\in(0,2]$ when its theoretical analogue $\rho$ is not defined. However here, $\rho_n$  does not follow the trend of the rank correlation coefficients.  Observe that for $t=0$, we have 
$$
P(X\leq x, Y\leq y)=P(X\leq x)\cdot 1\quad \mbox{and}\quad P(Y\leq y)=0.
$$
That way, the case  $t=0$ is not associated with independence, but $\rho_n$ wrongly indicates it. In general, the coefficient $\rho_n$ gives us  wrong information when $t\in(0,2]$ in spite of  the statistic $S^2_{\rho_{n}}$  takes  there on small values. The values of $\bar{\rho}_n$  are small for $t\in(0,2]$, because  for $t=0.1$ we have $cov_{sam}(X,Y)=-7.6*10\cdot^{30},\ S_x=0.9542$ and $S_y=2.5\cdot 10^{32}$, where $cov_{sam}$ is the sample covariance. The independence implies that the data is wildly dispersed  in both directions $x$ and $y$. In this example it happens only in one direction, along the variable $y$.

We also simulated the correlation coefficients  between correlation coefficients themselves. For $t=0.1$ we present below all  results for all  types of correlation coefficients:
$$
\bar{\rho}_{100}(\rho_{100}^{(i)}, \rho_{100,S}^{(j)})=-0.1265,\quad \bar{\rho}_{100}(\rho_{100}^{(i)}, \tau_{100}^{(j)})=-0.1334,\quad \bar{\rho}_{100}(\rho_{100}^{(i)}, r_{100}^{(j)})=-0.1360,
$$
$$
\bar{\rho}_{100}(\rho_{100,S}^{(i)}, \tau_{100}^{(j)})=0.9873,\quad \bar{\rho}_{100}(\rho_{100,S}^{(i)}, r_{100}^{(j)})=0.9660,\quad \bar{\rho}_{100}(\tau_{100}^{(i)}, r_{100}^{(j)})=0.9948,
$$
$$
\bar{\rho}_{100,S}(\rho_{100}^{(i)}, \rho_{100,S}^{(j)})=-0.1504,\quad \bar{\rho}_{100,S}(\rho_{100}^{(i)}, \tau_{100}^{(j)})=-0.1574,\quad \bar{\rho}_{100,S}(\rho_{100}^{(i)}, r_{100}^{(j)})=-0.1601,
$$
$$
\bar{\rho}_{100,S}(\rho_{100,S}^{(i)}, \tau_{100}^{(j)})=0.9842,\quad \bar{\rho}_{100,S}(\rho_{100,S}^{(i)}, r_{100}^{(j)})=0.9609,\quad \bar{\rho}_{100,S}(\tau_{100}^{(i)}, r_{100}^{(j)})=0.9929,
$$
$$
\bar{\tau}_{100}(\rho_{100}^{(i)}, \rho_{100,S}^{(j)})=-0.1021,\quad \bar{\tau}_{100}(\rho_{100}^{(i)}, \tau_{100}^{(j)})=-0.1068,\quad \bar{\tau}_{100}(\rho_{100}^{(i)}, r_{100}^{(j)})=-0.1088,
$$
$$
\bar{\tau}_{100}(\rho_{100,S}^{(i)}, \tau_{100}^{(j)})=0.9873,\quad \bar{\tau}_{100}(\rho_{100,S}^{(i)}, r_{100}^{(j)})=0.8393,\quad \bar{\tau}_{100}(\tau_{100}^{(i)}, r_{100}^{(j)})=0.9382.
$$
As we can see, the correlations between  the rank correlation coefficients are rather close to one. The  coefficient $\rho_n$ is even negatively correlated with the rank correlation coefficients. 

When $t$ increases the situation does not change much. For $t=1, 10, 100, 1000, 10000$,  we present the smallest and  largest  correlation coefficients:
$$
t=1,\quad \bar{\rho}_{100,S}(\rho_{100}^{(i)}, \rho_{100}^{(j)})=-0.0514,\quad \bar{\tau}_{100}(\tau_{100}^{(i)}, r_{100}^{(j)})=0.9947.
$$
$$
t=2,\quad \bar{\tau}_{100}(\rho_{100}^{(i)}, r_{100}^{(j)})=0.0634,\quad \bar{\rho}_{100}(\tau_{100}^{(i)}, r_{100}^{(j)})=0.9948.
$$
$$
t=10,\quad \bar{\rho}_{100,S}(\rho_{100}^{(i)}, r_{100,S}^{(j)})=0.4376,\quad \bar{\rho}_{100}(\tau_{100}^{(i)}, r_{100}^{(j)})=0.9948.
$$
$$
t=100,\quad \bar{\tau}_{100}(\rho_{100}^{(i)}, r_{100,S}^{(j)})=0.3870,\quad \bar{\rho}_{100}(\tau_{100}^{(i)}, r_{100}^{(j)})=0.9948.
$$
$$
t=1000,\quad \bar{\tau}_{100}(\rho_{100}^{(i)}, r_{100,S}^{(j)})=0.4037,\quad \bar{\rho}_{100}(\tau_{100}^{(i)}, r_{100}^{(j)})=0.9948.
$$
$$
t=10000,\quad \bar{\tau}_{100}(\rho_{100}^{(i)}, r_{100,S}^{(j)})=0.3903,\quad \bar{\rho}_{100}(\tau_{100}^{(i)}, r_{100}^{(j)})=0.9948.
$$

As was pointed out in Introduction, the coefficients  $\tau$ and $\rho_S$, and then, consequently, $r$ are  invariant under any increasing transformation. Let $1-\frac{1}{(1+y)^t}=v\in[0,1]$. Then
$$
F_1(x,v)=v-e^{-x}+(1-v)e^{-\frac{x}{1-v}}\quad (x\geq 0,\ 0\leq v\leq 1).
$$
For $F_1$, we have the same results for $\rho_S,\ \tau$ and $r$ as we had earlier for $F$, i.e., 
$$
\rho=-\frac{1}{\sqrt{3}}=-0.5774,\quad \rho_S=12\log 2-9=-0.6822,\quad \tau=-0.5,\quad r= 15/4-6\log 2=-0.4089.
$$
We  present the means  and variances of one thousand simulations of these correlation coefficients computed for samples $X_1,\ldots, X_{1000}$ and $Y_1,\ldots, Y_{1000}$ with $F_1$:
$$
\bar{\rho}_n=-0.5780,\quad \bar{\rho}_{n,S}=-0.6816,\quad \bar{\tau}_n=-0.4999,\quad \bar{r}_n=-0.4091,
$$
$$
S^2_{\rho_n}=3.2\cdot 10^{-4},\quad S^2_{\rho_{n,S}}=3.8\cdot 10^{-4},\quad S^2_{\tau_n}=2.7\cdot 10^{-4},\quad S^2_{r_n}=2.3\cdot 10^{-4}.
$$

\noindent {\bf Example~3.4}\
Let
$$
F(x,y)=\frac{x(y+1)(1+\alpha (1-x)y)}{2}\quad (0\leq x\leq 1,\ -1\leq y\leq 0)
$$
$$
F(x,y)=\frac{x}{2}+\frac{xy(1+\alpha (1-x)(1-y))}{2}\quad (0\leq x\leq 1,\ 0<y\leq 1),
$$
where $t \in[-1,1]$. If $t \not=0$, the variables $X$ and $Y$ are dependent, but  
$$
\tau =\tau_S=\rho=r =0.
$$
Between uncorrelatedness and independence lies semi-independence that is defined by the two equalities  $E(Y |X) = EY$ and $E(X|Y ) = EX$; see Jensen (1988). However here, $E(Y |X) = EY=0$, but 
$$
\frac{1}{2}-\frac{\alpha}{6}=E(X|Y ) \not= EX=\frac{1}{2}.
$$
The corresponding simulation results are presented in Tables~3.9--3.10 below.
\begin{table}[ht]
\begin{center}
{\em Table~3.9 Sample means of correlation coefficients for different $t$.}

\vspace*{1ex}
\begin{tabular}{|c|c|c|c|c|c|c|c|c|}
\hline
 &  $t=-1$ & $t=-0.5$ & $t=-0.1$ & $t=0$ & $t=0.1$ & $t=0.5$ & $t=1$\\
 \hline
 $\bar{\rho}_{n}$  & $-1.0\cdot 10^{-3}$ & $-2.1\cdot 10^{-4}$ & $8.6\cdot 10^{-4}$ & $-2.0\cdot 10^{-3}$ & $3.5\cdot 10^{-4}$ & $-8.5\cdot 10^{-4}$ & $1.6\cdot 10^{-3}$ \\
\hline
$\bar{\rho}_{n,S}$ & $-1.3\cdot 10^{-3}$ & $-2.1\cdot 10^{-3}$ & $9.2\cdot 10^{-4}$ & $-2.0\cdot 10^{-3}$ & $3.0\cdot 10^{-4}$ & $-7.6\cdot 10^{-4}$ & $1.4\cdot 10^{-3}$ \\
\hline
$\bar{\tau}_{n}$   & $-9.2\cdot 10^{-4}$ & $-1.4\cdot 10^{-3}$ & $6.1\cdot 10^{-4}$ & $-1.8\cdot 10^{-3}$ & $2.2\cdot 10^{-4}$ & $-4.9\cdot 10^{-4}$ & $9.9\cdot 10^{-4}$ \\
\hline
 $\bar{r}_{n}$     & $-7.2\cdot 10^{-4}$ & $-1.0\cdot 10^{-3}$ & $4.6\cdot 10^{-4}$ & $-1.0\cdot 10^{-3}$ & $1.7\cdot 10^{-4}$ & $-3.6\cdot 10^{-4}$ & $7.7\cdot 10^{-4}$ \\
\hline
\end{tabular}
\end{center}
\end{table}

\begin{table}[ht]
\begin{center}
{\em Table~3.10 Sample variances of correlation coefficients for different $t$.}

\vspace*{1ex}
\begin{tabular}{|c|c|c|c|c|c|c|c|c|}
\hline
 &  $t=-1$ & $t=-0.5$ & $t=-0.1$ & $t=0$ &  $t=0.1$ & $t=0.5$ & $t=1$\\
 \hline
 $S^2_{\rho_{n}}$  & $9.7\cdot 10^{-4}$ & $9.4\cdot 10^{-4}$ & $9.9\cdot 10^{-4}$ & $9.4\cdot 10^{-4}$ & $9.7\cdot 10^{-4}$ & $1.1\cdot 10^{-3}$ & $9.9\cdot 10^{-4}$ \\
\hline
$S^2_{\rho_{n,S}}$ & $1.0\cdot 10^{-3}$ & $9.5\cdot 10^{-4}$ & $9.9\cdot 10^{-4}$ & $9.4\cdot 10^{-4}$ & $9.7\cdot 10^{-4}$ & $1.1\cdot 10^{-3}$ & $1.0\cdot 10^{-3}$ \\
\hline
$S^2_{\tau_n}$     & $4.6\cdot 10^{-4}$ & $4.2\cdot 10^{-4}$ & $4.4\cdot 10^{-4}$ & $4.2\cdot 10^{-4}$ & $4.3\cdot 10^{-4}$ & $4.7\cdot 10^{-4}$ & $4.7\cdot 10^{-4}$ \\
\hline
  $S^2_{r_n}$      & $2.6\cdot 10^{-4}$ & $2.4\cdot 10^{-4}$ & $2.5\cdot 10^{-4}$ & $2.4\cdot 10^{-4}$ & $2.4\cdot 10^{-4}$ & $2.6\cdot 10^{-4}$ & $2.7\cdot 10^{-4}$ \\
\hline
\end{tabular}
\end{center}
\end{table} 
The coefficient  $r_n$  has the smallest sample variance and can again be viewed as a good rank correlation coefficient. It is interesting that the generation of correlation coefficients for $t=0.5$ gave us  negative numbers.

\section*{Conclusion}\ In the present paper, we have discussed the Pearson $\rho$,  Spearman $\rho_S$,  Kendall $\tau$ correlation coefficients  and have introduced a new correlation coefficient $r$.  We have proposed to consider  $\tilde{\rho}$, such that  $\rho_n$ approximates $\tilde{\rho}$ in simulation experiments for large $n$,
as a new extension of $\rho$. 

We have conducted  simulation experiments and studied the behavior of the above correlation coefficients. We have come to the conclusion that the behavior of $\rho_n$ can be very different from the behavior of the rank correlation coefficients $\rho_{n,S},\ \tau_n,\ r_n$. This follows from the definition of $\rho_n$ and $\rho_{n,S},\ \tau_n\ r_n$. The coefficient $\rho_n$ measures the association between variables $X_i$ and $Y_j$, whereas  the coefficients $\rho_{n,S},\ \tau_n$ and $r_n$ measure the association between their ranks. In this paper we have raised the question:  which  correlation coefficient measures the dependence rate better? We think that it is fair to compare  only    $\tau_n,\ \rho_{n,S}$ and $r_n$, since these three coefficients   are the ranks correlation coefficients. We  have come  to the conclusion that $r$ looks a little better than $\rho_S$ and $\tau$. The coefficient $\rho$, from one hand, and the coefficients $\rho_{S},\ \tau,\ r$, from the other hand,  measure in their own peculiar ways, and cannot be compared generally. We recommend to use $\rho$ when the association between $X$ and $Y$ is close to $\pm 1$ and to use $\rho_{S},\ \tau,\ r$, otherwise. 

\section*{References}
\begin{description} 
{\small

\item Bairamov, I., Stepanov, A. (2010).\ Numbers of near-maxima for the bivariate case,  {\it Statistics  $\&$ Probability Letters}, {\bf 80}, 196--205.

\item Balakrishnan, N., Lai, C.D. (2009).\ {\it Continuous Bivariate Distributions}, Second edition, Springer.
 
\item Balakrishnan, N., Stepanov, A. (2015).\ Limit results for concomitants of order statistics, {\it Metrika}, {\bf 78}, 385--397.

\item  Bhattacharya, B.B. (1974).\ Convergence of sample paths of normalized sums of induced order statistics, {\it Ann. Statist. }, {\bf 2}, 1034--1039.

\item Bhattacharya, B.B. (1984).\ Induced order statistics: theory and applications, In {\it Handbook of Statistics 4}, Ed.
Krishnaiah, P. R., Sen, P. K.,  North Holland, Amsterdam, 383--403. 

\item Chu, S.J.,  Huang, W.J., Chen, H. (1999).\ A study of asymptotic distributions of concomitants of certain order statistics, {\it Statist. Sinica}, {\bf 9}, 811--830.

\item  Cifarelli, D. M., Conti, P. L., Regazzini, E. (1996).\ On the asymptotic distribution of a general measure of monotone dependence, {\it The Annals of Statistics}, {\bf 25}, 	1386--1399.

\item Daniels, H. E. (1950).\ Rank correlation and population models, {\it Journal of the Royal Statistical Society}, Ser. B, {\bf 12} (2),
171--191.

\item  David, H.A. (1994).\ Concomitants of Extreme Order Statistics, In {\it Extreme Value Theory and Applications},
Proceedings of the Conference on Extreme Value Theory and Applications, {\bf 1}, Ed. Galambos, J., Lechner, J., Simiu, E., Kluwer Academic Publishers, Boston  211--224.

\item David, H.A., Galambos, J. (1974).\ The asymptotic theory of concomitants of order statistics, {\it J. Appl. Probab.}, {\bf 11}, 762--770.

\item David, H.A., Nagaraja, H.N. (2003).\  {\it Order Statistics}, Third  edition, John Wiley \& Sons, NY.

\item Durbin J., Stuart,  A. (1951).\ Inversions and rank correlation coefficients, {\it Journal of the Royal Statistical Society},  Ser. B {\bf 13} (2), 303--309.

\item Egorov, V. A., Nevzorov, V. B. (1984).\ Rate of convergence to the Normal law of sums of induced order statistics,
{\it Journal of Soviet Mathematics} (New York), {\bf 25}, 1139--1146.

\item Goel, P. K., Hall, P. (1994).\ On the average difference between concomitants and order statistics,
{\it  Ann. Probab.}, {\bf 22}, 126--144.

\item Fisher, R. A. (1921).\ On the ''probable error" of a coefficient of correlation deduced from a small sample, Metron, {\bf 1},  3--32.

\item Jensen, D.R. (1988).\  Semi-independence. In: Encyclopedia of Statistical Sciences, Volume 8, S. Kotz and N.L. Johnson (eds.),  358–-359. John Wiley and Sons, New York?

\item Kendall, M. G. (1970).\ {\it Rank Correlation Methods}, London, Griffin.

\item Omey, E., Gulk, S.Van (2008).\ Central limit theorems for variances and correlation coefficients, {\it Preprint}, doi: 10.13140/RG.2.2.21663.66727.

\item Stepanov, A.V. (2024).\ On Kendall's correlation coefficient, {\it Vestnik St. Petersburg University, Mathematics}, accepted.

\item  Xu, W., Hou,  Y.,   Hung, Y. S., Zou, Y. (2009). Comparison of Spearman's rho and Kendall's
tau in Normal and Contaminated Normal Models, {\it arXiv:1011.2009v1 [cs.IT]}.

}

\end{description}

\end{document}